# LYAPUNOV THEOREMS FOR SYSTEMS DESCRIBED BY RETARDED FUNCTIONAL DIFFERENTIAL EQUATIONS

**Iasson Karafyllis**
Division of Mathematics, Dept. of Economics, University of Athens
Email: ikarafil@econ.uoa.gr

**Abstract**
Lyapunov-like characterizations for non-uniform in time and uniform robust global asymptotic stability of uncertain systems described by retarded functional differential equations are provided.

**Keywords:** Lyapunov functionals, time-delay systems, robust global asymptotic stability.

## 1. Introduction

In this paper we provide Lyapunov characterizations for non-uniform in time and uniform Robust Global Asymptotic Stability (RGAS) for systems described by time-varying Retarded Functional Differential Equations (RFDEs). The notion of non-uniform in time RGAS is introduced in [10] for continuous time finite-dimensional systems and in [12] for a wide class of systems including discrete-time systems and systems described by RFDEs. This notion has been proved to be fruitful for the solution of several problems in Mathematical Control Theory (see [10,11]) and is a special case of the concept of stability with respect to two measures introduced in [17]. The notion of uniform RGAS that we adopt in this paper is an extension of the corresponding notion for finite-dimensional continuous-time uncertain systems (see [3,19]).

The motivation for the extension of non-uniform in time and uniform RGAS to uncertain systems described by RFDEs is strong, since such models are used frequently for the description of engineering systems (see [8]). It should be emphasized that in many cases where hybrid open-loop/feedback stabilizing control laws are proposed for finite-dimensional continuous-time systems, the closed-loop system is actually a system described by time-varying RFDEs (infinite-dimensional). For example, in [24] analytic driftless control systems of the following form are considered:

$$\dot{x}(t) = f(x(t), u(t)) := \sum_{i=1}^{m} f_i(x(t)) u_i(t)$$

$$x(t) \in \Re^n, \ u(t) := (u_1(t), ..., u_m(t))' \in \Re^m$$

The authors in [24] provide strategies for the construction of control laws of the form $u(t) = k(t, x(t), x(lr))$ for $t \in [lr, (l+1)r)$, where $l$ is a non-negative integer and $r > 0$ denotes the updating time-period of the control. Notice that the closed-loop system is described by the equations:

$$\dot{x}(t) = f\left(x(t), k\left(t, x(t), x\left(\left[\frac{t}{r}\right]r\right)\right)\right)$$

$$x(t) \in \Re^n, \ t \geq 0$$

where $\left[\frac{t}{r}\right]$ denotes the integer part of $\frac{t}{r}$. The above system is actually a system described by RFDEs and is time-varying even if $k$ is independent of time. The same comments apply for the synchronous controller switching strategies proposed in [26]. The possibility of switching control laws using distributed delays was exploited in [22]. However, it should be emphasized that hybrid systems with asynchronous switching or state space depending switching rules (see for example [4,26]) in general cannot be described by RFDEs with Lipschtiz continuous right-hand sides.

Lyapunov functions and functionals play an important role to synthesis and design in control theory and several important results have been established concerning Lyapunov-like descriptions of *uniform global asymptotic stability* (UGAS) (see [6,7,8,16] and the references therein). Our goal is to establish Lyapunov characterizations for the concepts of *non-uniform in time and uniform robust global asymptotic stability* (RGAS) analogous to the corresponding characterizations given in [3,19] for continuous-time finite-dimensional uncertain systems, which overcome the limitations imposed by previous works. Particularly, our Lyapunov characterizations apply

- to systems with disturbances that take values in a (not necessarily compact) given set
- to systems described by RFDEs with right-hand sides which are not necessarily bounded with respect to time

The difficulties in verifying stability properties for the time-varying case with right-hand sides, which are not necessarily bounded with respect to time, are explained in [6,7]. Our motivation to provide Lyapunov characterizations that cover the above case is strong and is explained below:

• as in the finite-dimensional case, where the Lyapunov characterization for systems with disturbances given in [19] led to Lyapunov characterizations of the Input-to-State Stability (ISS, see [27]) property, we believe that our results will similarly lead to Lyapunov characterizations of the ISS property for systems described by RFDEs
• as in the finite-dimensional case, where the Lyapunov characterization given in [19] led to Lyapunov-like conditions for the robust stabilization of finite-dimensional control systems by means of continuous static feedback, we believe that our results can be used for the expression of Lyapunov-like conditions for the robust stabilization of systems described by RFDEs by means of continuous distributed delay static feedback
• we believe that our results can be used for the study of the robustness properties with respect to modeling, actuator and measurement errors of the closed-loop system for finite-dimensional control systems under hybrid feedback control strategies or continuous distributed delay static feedback (as explained previously)
• as in the finite-dimensional case, where it was shown that certain control systems cannot be uniformly stabilized by means of time-invariant feedback laws but can be stabilized non-uniformly in time by means of time-varying static feedback (see [10,11]), we believe that our results can be used for exploitation of the stabilizing capabilities of time-varying distributed delay feedback and this clearly motivates the study of **non-uniform in time RGAS** and its Lyapunov characterizations

It should be emphasized that the studies described above **cannot** be performed using the existing Lyapunov characterizations of stability for systems described by RFDEs.

Viability issues for systems described by functional differential inclusions (and thus uncertain systems described by RFDEs) were considered in [2]. Lyapunov-like conditions that guarantee stability with respect to part of the variables for systems described by RFDEs are provided in [30] (Chapter 6). We note that Lyapunov functionals for linear time-delay systems were constructed in [5,15,23]. Stability conditions are given in [6,7,14,16] for time-varying time-delay systems and in [18,31] for time-varying integrodifferential systems. Recently in many works the problem of feedback stabilization of systems described by RFDEs was studied (see for instance [9,20,21,25]). It should be emphasized that the literature concerning issues of stability and stabilization of linear time-delay systems is vast and the previous references are only given as pointers. Note also that in the present paper we are not concerned with stability conditions given by Razumikhin functions, since such conditions resemble "small-gain" conditions with Lyapunov-like characteristics (see [29]).

In the present work we provide Lyapunov-like conditions which demand the infinitesimal decrease property to hold only on subsets of the state space which contain the solutions of the system (i.e., the infinitesimal decrease property holds only after some time) along with an additional property that guarantees forward completeness on the critical time interval when the infinitesimal decrease property does not hold, namely the property $\dot{V} \leq V$, where $V$ denotes the Lyapunov functional and $\dot{V}$ the time derivative of the Lyapunov functional evaluated along the solutions of the system (Theorems 2.9 and 2.10). This property was shown to be necessary and sufficient for forward completeness of continuous-time finite dimensional systems in [1]. Moreover, the "weaker" property that we demand in the present paper is utilized for the construction of Lyapunov functionals for time-delay systems (see Examples 2.12 and 2.13).

**Notations** Throughout this paper we adopt the following notations:
∗ For a vector $x \in \Re^n$ we denote by $|x|$ its usual Euclidean norm and by $x'$ its transpose. For $x \in C^0([-r,0];\Re^n)$ we define $\|x\|_r := \max_{\theta \in [-r,0]} |x(\theta)|$.

∗ We denote by $[R]$ the integer part of the real number $R$, i.e., the greatest integer, which is less than or equal to $R$.

∗ By $C^0(A;\Omega)$, we denote the class of continuous functions on $A$, which take values in $\Omega$.

∗ $\mathcal{E}$ denotes the class of non-negative $C^0$ functions $\mu: \Re^+ \to \Re^+$, for which it holds: $\int_0^{+\infty} \mu(t)dt < +\infty$ and $\lim_{t \to +\infty} \mu(t) = 0$.

* $Z^+$ denotes the set of positive integers.
* We denote by $K^+$ the class of positive $C^\infty$ functions defined on $\Re^+$. We say that a function $\rho : \Re^+ \to \Re^+$ is positive definite if $\rho(0) = 0$ and $\rho(s) > 0$ for all $s > 0$. By $K$ we denote the set of positive definite, increasing and continuous functions. We say that a positive definite, increasing and continuous function $\rho : \Re^+ \to \Re^+$ is of class $K_\infty$ if $\lim_{s \to +\infty} \rho(s) = +\infty$. By $KL$ we denote the set of all continuous functions $\sigma = \sigma(s,t) : \Re^+ \times \Re^+ \to \Re^+$ with the properties: (i) for each $t \geq 0$ the mapping $\sigma(\cdot, t)$ is of class $K$; (ii) for each $s \geq 0$, the mapping $\sigma(s, \cdot)$ is non-increasing with $\lim_{t \to +\infty} \sigma(s,t) = 0$.

## 2. Robust Global Asymptotic Stability (RGAS) for Systems Described by RFDEs

This section of the paper is structured as follows. In the first subsection (Subsection 2.I) we provide all the background material that the reader needs to know in order to understand the results of the present paper. In the next subsection (Subsection 2.II) the reader is introduced to the stability notions used in the present paper as well as to some preliminary results already presented in [12]. Subsection 2.III is devoted to the presentation of some tools of non-smooth analysis that are used for the expression of the infinitesimal decrease property of the Lyapunov functionals, while Subsection 2.IV is devoted to the statement of our main results and the presentation of examples.

### 2.I. Background Material on Systems Described by RFDEs

Let $x : [a-r, b) \to \Re^n$ with $b > a \geq 0$ and $r \geq 0$. We define for $t \in [a,b)$

$$T_r(t)x := x(t+\theta) \, ; \, \theta \in [-r, 0] \tag{2.1}$$

Let $D \subseteq \Re^l$ a non-empty set and $M_D$ the class of all right-continuous mappings $d : \Re^+ \to D$, with the following property:

"there exists a countable set $A_d \subset \Re^+$ which is either finite or $A_d = \{t_k^d \, ; \, k = 1,...,\infty\}$ with $t_{k+1}^d > t_k^d > 0$ for all $k = 1,2,...$ and $\lim t_k^d = +\infty$, such that the mapping $t \in \Re^+ \setminus A_d \to d(t) \in D$ is continuous"

We denote by $x(t)$ with $t \geq t_0$ the unique solution of the initial-value problem:

$$\begin{aligned} \dot{x}(t) &= f(t, T_r(t)x, d(t)) \, , \, t \geq t_0 \\ x(t) &\in \Re^n \, , \, d(\cdot) \in M_D \end{aligned} \tag{2.2}$$

with initial condition $T_r(t_0)x = x_0 \in C^0([-r,0]; \Re^n)$, where $r \geq 0$ is a constant and the mapping $f : \Re^+ \times C^0([-r,0]; \Re^n) \times D \to \Re^n$ with $f(t,0,d) = 0$ for all $(t,d) \in \Re^+ \times D$ satisfies the following hypotheses:

**(H1)** The mapping $(x,d) \to f(t,x,d)$ is continuous for each fixed $t \geq 0$ and such that for every bounded $I \subseteq \Re^+$ and for every bounded $S \subset C^0([-r,0]; \Re^n)$, there exists a constant $L \geq 0$ such that:

$$(x(0) - y(0))'(f(t,x,d) - f(t,y,d)) \leq L \max_{\tau \in [-r,0]} |x(\tau) - y(\tau)|^2 = L\|x-y\|_r^2$$

$$\forall t \in I \, , \, \forall (x,y) \in S \times S \, , \, \forall d \in D$$

This assumption is equivalent to the existence of a continuous function $L : \Re^+ \times \Re^+ \to \Re^+$ such that for each fixed $t \geq 0$ the mappings $L(t, \cdot)$ and $L(\cdot, t)$ are non-decreasing, with the following property:

$$\left(x(0)-y(0)\right)'\left(f(t,x,d)-f(t,y,d)\right) \leq L(t,\|x\|_r+\|y\|_r)\|x-y\|_r^2 \quad (2.3)$$
$$\forall (t,x,y,d) \in \Re^+ \times C^0([-r,0];\Re^n) \times C^0([-r,0];\Re^n) \times D$$

**(H2)** For every bounded $\Omega \subset \Re^+ \times C^0([-r,0];\Re^n)$ the image set $f(\Omega \times D) \subset \Re^n$ is bounded.

**(H3)** There exists a countable set $A \subset \Re^+$, which is either finite or $A = \{t_k\; ;\; k=1,...,\infty\}$ with $t_{k+1} > t_k > 0$ for all $k = 1,2,...$ and $\lim t_k = +\infty$, such that mapping $(t,x,d) \in (\Re^+ \setminus A) \times C^0([-r,0];\Re^n) \times D \to f(t,x,d)$ is continuous. Moreover, for each fixed $(t_0,x,d) \in \Re^+ \times C^0([-r,0];\Re^n) \times D$, we have $\lim_{t \to t_0^+} f(t,x,d) = f(t_0,x,d)$.

**(H4)** For every $\varepsilon > 0$, $t \in \Re^+$, there exists $\delta := \delta(\varepsilon,t) > 0$ such that $\sup\left\{|f(\tau,x,d)|\; ;\; \tau \in \Re^+, d \in D, |\tau-t|+\|x\|_r < \delta \right\} < \varepsilon$.

It is clear that for every $d \in M_D$ the composite map $f(t,x,d(t))$ satisfies the Caratheodory condition on $\Re^+ \times C^0([-r,0];\Re^n)$ and consequently, by virtue of Theorem 2.1 in [8] (and its extension given in paragraph 2.6 of the same book), for every $(t_0,x_0,d) \in \Re^+ \times C^0([-r,0];\Re^n) \times M_D$ there exists $h > 0$ and at least one continuous function $x : [t_0 - r, t_0 + h] \to \Re^n$, which is absolutely continuous on $[t_0, t_0 + h]$ with $T_r(t_0)x = x_0$ and $\dot{x}(t) = f(t, T_r(t)x, d(t))$ almost everywhere on $[t_0, t_0+h]$. Let $x:[t_0-r,t_0+h] \to \Re^n$ and $y:[t_0-r,t_0+h] \to \Re^n$ two solutions of (2.2) with initial conditions $T_r(t_0)x = x_0$ and $T_r(t_0)y = y_0$ and corresponding to the same $d \in M_D$. Evaluating the derivative of the absolutely continuous map $z(t) = |x(t)-y(t)|^2$ on $[t_0,t_0+h]$ in conjunction with hypothesis (H1) above, we obtain the integral inequality:

$$|x(t)-y(t)|^2 \leq |x(t_0)-y(t_0)|^2 + 2\int_{t_0}^t \tilde{L}\|T_r(\tau)x - T_r(\tau)y\|_r^2 d\tau, \quad \forall t \in [t_0,t_0+h]$$

where $\tilde{L} := L(t,a(x,y))$, $L(\cdot)$ is the function involved in (2.3) and $a(x,y) := \sup_{t \in [t_0-r,t_0+h]} |x(t)| + \sup_{t \in [t_0-r,t_0+h]} |y(t)|$. Consequently, we obtain:

$$\|T_r(t)(x-y)\|_r^2 \leq \|x_0-y_0\|_r^2 + 2\int_{t_0}^t \tilde{L}\|T_r(\tau)(x-y)\|_r^2 d\tau, \quad \forall t \in [t_0,t_0+h]$$

and immediate application of the Gronwall-Bellman inequality gives:

$$\|T_r(t)(x-y)\|_r \leq \|x_0-y_0\|_r \exp\left(\tilde{L}(t-t_0)\right), \quad \forall t \in [t_0,t_0+h] \quad (2.4)$$

Thus, we conclude that under hypotheses (H1-4), for every $(t_0,x_0,d) \in \Re^+ \times C^0([-r,0];\Re^n) \times M_D$ there exists $h > 0$ and exactly one continuous function $x:[t_0-r,t_0+h] \to \Re^n$, which is absolutely continuous on $[t_0,t_0+h]$ with $T_r(t_0)x = x_0$ and $\dot{x}(t) = f(t,T_r(t)x,d(t))$ almost everywhere on $[t_0,t_0+h]$. We denote by $\phi(t,t_0,x_0;d) := T_r(t)x$ and $\phi(\theta,t,t_0,x_0;d) := x(t+\theta)$ for $\theta \in [-r,0]$. Clearly, we have $\phi(t,t_0,x_0;d) = \phi(\theta,t,t_0,x_0;d)\; ;\; \theta \in [-r,0]$. Clearly, the unique solution of (2.2) satisfies for all $(t_0,x_0,y_0,d) \in \Re^+ \times C^0([-r,0];\Re^n) \times C^0([-r,0];\Re^n) \times M_D$ and for all $t \geq t_0$ so that $\phi(t,t_0,x_0;d)$ and $\phi(t,t_0,y_0;d)$ are both defined:

$$\|\phi(t,t_0,x_0;d) - \phi(t,t_0,y_0;d)\|_r \leq \|x_0-y_0\|_r \exp\left(\tilde{L}(t-t_0)\right) \quad (2.5)$$

where $\widetilde{L} := L\left(t, \sup_{t\in[t_0,t]} \|\phi(t,t_0,x_0;d)\|_r + \sup_{t\in[t_0,t]} \|\phi(t,t_0,y_0;d)\|_r\right)$ and $L(\cdot)$ is the function involved in (2.3).

Using hypothesis (H2) above and Theorem 3.2 in [8], we conclude that for every $(t_0,x_0,d) \in \Re^+ \times C^0([-r,0];\Re^n) \times M_D$ there exists $t_{\max} \in (t_0,+\infty]$, such that the unique solution $x(t)$ of (2.2) is defined on $[t_0-r,t_{\max})$ and cannot be further continued. Moreover, if $t_{\max} < +\infty$ then we must necessarily have $\limsup_{t\to t_{\max}^-}|x(t)| = +\infty$. In addition, it is clear that for every $d \in M_D$ the composite map $f(t,x,d(t))$ is continuous on $(\Re^+ \setminus (A \cup A_d)) \times C^0([-r,0];\Re^n)$. Applying repeatedly Theorem 2.1 in [8] on each one of the intervals contained in $[t_0,t_{\max}) \setminus (A \cup A_d)$, we conclude that the solution satisfies $\dot{x}(t) = f(t,T_r(t)x,d(t))$ for all $t \in [t_0,t_{\max}) \setminus (A \cup A_d)$. Since, the composite map $t \to f(t,x,d(t))$ is right-continuous on $\Re^+$, by virtue of the mean value theorem, it follows that $\lim_{h\to 0^+} \frac{x(t+h)-x(t)}{h} = f(t,T_r(t)x,d(t))$ for all $t \in [t_0,t_{\max})$.

**Remark 2.1:**

(a) When $r = 0$ we identify the space $C^0([-r,0];\Re^n)$ with the finite-dimensional space $\Re^n$ and we obtain the familiar finite-dimensional continuous-time case. Consequently, all the following results hold also for finite-dimensional continuous-time systems.

(b) A major difference between the case of uncertain finite-dimensional continuous-time systems considered in [19] and the case of uncertain systems described by RFDEs is the nature of the class of allowed inputs $M_D$. This happens because there is a fundamental difference between the two cases: in the finite-dimensional case the map describing the evolution of the state is absolutely continuous with respect to time while in the infinite-dimensional case the map describing the evolution of the state is (simply) continuous with respect to time (see Lemma 2.1 in [8] and notice that the state for the infinite-dimensional case is $T_r(t)x \in C^0([-r,0];\Re^n)$). This fact has an important consequence: Lyapunov functionals evaluated on the solutions of system (2.2) will be (simply) continuous and not absolutely continuous maps with respect to time and in order to guarantee their monotonicity, we must require that an appropriate decrease condition holds for all times (and not almost everywhere, see the discussion in [3], Chapter 6). Thus we cannot allow $M_D$ contain arbitrary measurable mappings.

(c) In all the following results we assume that the inputs belong to the class $M_D$. It is clear that the same conclusions hold for inputs $d: \Re^+ \to D$, for which there exists $d' \in M_D$ such that $d(t) = d'(t)$ almost everywhere.

(d) As already pointed out in the Introduction, it should be emphasized that model (2.2) under hypotheses (H1-4) can represent finite-dimensional control systems $\dot{x}(t) = f(t,x(t),u(t)), (t,x(t),u(t)) \in \Re^+ \times \Re^n \times \Re^m$ under hybrid feedback control with synchronous switching of the form $u(t) = k(t,x(t),x(t_i))$, for $t_i \le t < t_{i+1}$, where $\pi = \{t_i, i=0,1,2,...\}$ is a partition of $\Re^+$ of diameter $r > 0$ (i.e., $\pi = \{t_i, i=0,1,2,...\}$ is an increasing sequence with $t_0 = 0$, $\lim t_i = +\infty$ and $t_{i+1} - t_i \le r$ for all $i = 0,1,2,...$), when the vector fields $f(t,x,u)$ and $k(t,x,x')$ are continuous and locally Lipschitz with respect to $(x,u)$ and $(x,x')$, respectively, with $f(t,0,k(t,0,0)) = 0$ for all $t \ge 0$. Particularly, if we define the function $p(t) = \max\{\tau; \tau \in \pi, \tau \le t\}$ and the mapping $(t,x) \in \Re^+ \times C^0([-r,0];\Re^n) \to \widetilde{f}(t,x) \in \Re^n$, where $\widetilde{f}(t,x) = f(t,x(0),k(t,x(0),x(p(t)-t)))$, then $\widetilde{f}$ satisfies hypotheses (H1-4) and moreover the closed-loop system is described by the RFDEs $\dot{x}(t) = \widetilde{f}(t,T_r(t)x)$. Notice that if $f$ and $k$ are independent of $t \ge 0$ (time-invariant vector fields), the mapping $\widetilde{f}(t,x) = f(x(0),k(x(0),x(p(t)-t)))$ is time-varying. Moreover, for the special case $t_i = ir$, $i = 0,1,2,...$, then $p(t) = \left[\frac{t}{r}\right]r$ and the mapping $\widetilde{f}(t,x) = f\left(x(0),k\left(x(0),x\left(\left[\frac{t}{r}\right]r - t\right)\right)\right)$ is periodic with respect to $t \ge 0$ with period $r > 0$.

## 2.II. Background Material on Non-Uniform in Time and Uniform RGAS for Systems Described by RFDEs

Since $f(t,0,d) = 0$ for all $(t,d) \in \Re^+ \times D$, it follows that $\phi(t,t_0,0;d) = 0 \in C^0([-r,0];\Re^n)$ for all $(t_0,d) \in \Re^+ \times M_D$ and $t \geq t_0$. Furthermore, for every $\varepsilon > 0$, $T, h \geq 0$ there exists $\delta := \delta(\varepsilon,T,h) > 0$ such that:

$$\|x\|_r < \delta \quad \Rightarrow \quad \sup\{\|\phi(\tau,t_0,x;d)\|_r \ ; \ d \in M_D \ , \tau \in [t_0, t_0+h] \ , t_0 \in [0,T]\} < \varepsilon$$

Thus $0 \in C^0([-r,0];\Re^n)$ is an equilibrium point for (2.2) in the sense described in [12]. The following definition of non-uniform in time RGAS coincides with the definition of non-uniform in time RGAS given in [12], for a wide class of systems that include systems of RFDEs studied in this paper.

**Definition 2.2:** *We say that* $0 \in C^0([-r,0];\Re^n)$ *is **non-uniformly in time Robustly Globally Asymptotically Stable (RGAS)** for system (2.2) if the following properties hold:*

**P1** $0 \in C^0([-r,0];\Re^n)$ is **Robustly Lagrange Stable**, i.e., *for every* $s \geq 0$, $T \geq 0$, *it holds that*

$$\sup\{\|\phi(t,t_0,x_0;d)\|_r \ ; t \in [t_0,+\infty) \ , \|x_0\|_r \leq s \ , t_0 \in [0,T] \ , d \in M_D \ \} < +\infty$$
**(Robust Lagrange Stability)**

**P2** $0 \in C^0([-r,0];\Re^n)$ is **Robustly Lyapunov Stable**, i.e., *for every* $\varepsilon > 0$ *and* $T \geq 0$ *there exists a* $\delta := \delta(\varepsilon,T) > 0$ *such that:*

$$\|x_0\|_r \leq \delta \ , t_0 \in [0,T] \Rightarrow \|\phi(t,t_0,x_0;,d)\|_r \leq \varepsilon \ , \forall t \in [t_0,+\infty), \ \forall d \in M_D$$
**(Robust Lyapunov Stability)**

**P3** $0 \in C^0([-r,0];\Re^n)$ satisfies the **Robust Attractivity Property**, i.e. *for every* $\varepsilon > 0$, $T \geq 0$ *and* $R \geq 0$, *there exists a* $\tau := \tau(\varepsilon,T,R) \geq 0$, *such that:*

$$\|x_0\|_r \leq R \ , t_0 \in [0,T] \Rightarrow \|\phi(t,t_0,x_0;d)\|_r \leq \varepsilon \ , \forall t \in [t_0+\tau,+\infty), \ \forall d \in M_D$$

The two following lemmas are given in [12] (as Lemma 3.3 and Lemma 3.4, respectively) for a wide class of systems that include systems of RFDEs under hypotheses (H1-4). They provide essential characterizations of the notion of non-uniform in time RGAS.

**Lemma 2.3:** *Suppose that (2.2) is Robustly Forward Complete, i.e., for every* $s \geq 0$, $T \geq 0$, *it holds that*

$$\sup\{\|\phi(t_0+h,t_0,x_0;d)\|_r \ ; \ h \in [0,T] \ , \|x_0\|_r \leq s \ , t_0 \in [0,T] \ , d \in M_D \ \} < +\infty$$

*and that* $0 \in C^0([-r,0];\Re^n)$ *satisfies the Robust Attractivity Property (property P3 of Definition 2.2) for system (2.2). Then* $0 \in C^0([-r,0];\Re^n)$ *is non-uniformly in time RGAS for system (2.2).*

**Lemma 2.4:** $0 \in C^0([-r,0];\Re^n)$ *is non-uniformly in time RGAS for system (2.2) if and only if there exist functions* $\sigma \in KL$, $\beta \in K^+$ *such that the following estimate holds for all* $(t_0, x_0, d) \in \Re^+ \times C^0([-r,0];\Re^n) \times M_D$ *and* $t \in [t_0,+\infty)$:

$$\|\phi(t,t_0,x_0;d)\|_r \ \leq \sigma\big(\beta(t_0)\|x_0\|_r \ , t-t_0 \big) \tag{2.6}$$

Finally, we also provide the definition of uniform RGAS, in terms of $KL$ functions, which is completely analogous to the finite-dimensional case (see [3,19]). It is clear that such a definition is equivalent to a $\delta - \varepsilon$ definition (analogous to Definition 2.2).

**Definition 2.5:** *We say that* $0 \in C^0([-r,0]; \Re^n)$ *is Uniformly Robustly Globally Asymptotically Stable (URGAS) for system (2.2) if and only if there exist a function* $\sigma \in KL$ *such that the following estimate holds for all* $(t_0, x_0, d) \in \Re^+ \times C^0([-r,0]; \Re^n) \times M_D$ *and* $t \in [t_0, +\infty)$:

$$\|\phi(t, t_0, x_0; d)\|_r \leq \sigma\left(\|x_0\|_r, t - t_0\right) \tag{2.7}$$

The following corollary must be compared to Lemma 1.1, page 131 in [8]. It shows that for periodic systems of RFDEs non-uniform in time RGAS is equivalent to URGAS. We say that (2.2) is $T-periodic$ if there exists $T > 0$ such that $f(t+T, x, d) = f(t, x, d)$ for all $(t, x, d) \in \Re^+ \times C^0([-r,0]; \Re^n) \times D$.

**Corollary 2.6:** *Suppose that* $0 \in C^0([-r,0]; \Re^n)$ *is non-uniformly in time RGAS for system (2.2) and that (2.2) is* $T-periodic$. *Then* $0 \in C^0([-r,0]; \Re^n)$ *is URGAS for system (2.2).*

**Proof** The proof is based on the following observation: if (2.2) is $T-periodic$ then for all $(t_0, x_0, d) \in \Re^+ \times C^0([-r,0]; \Re^n) \times M_D$ it holds that $\phi(t, t_0, x_0; d) = \phi\left(t - \left[\frac{t_0}{T}\right]T, t_0 - \left[\frac{t_0}{T}\right]T, x_0; P(t_0)d\right)$, where $[t_0 / T]$ denotes the integer part of $t_0 / T$ and $P(t_0)d \in M_D$ is defined by:

$$(P(t_0)d)(t) := d\left(t + \left[\frac{t_0}{T}\right]T\right), \quad \forall t \geq 0$$

Since $0 \in C^0([-r,0]; \Re^n)$ is non-uniformly in time RGAS for system (2.2), there exist functions $\sigma \in KL$, $\beta \in K^+$ such that (2.6) holds for all $(t_0, x_0, d) \in \Re^+ \times C^0([-r,0]; \Re^n) \times M_D$ and $t \in [t_0, +\infty)$. Consequently, it follows that the following estimate holds for all $(t_0, x_0, d) \in \Re^+ \times C^0([-r,0]; \Re^n) \times M_D$ and $t \in [t_0, +\infty)$:

$$\|\phi(t, t_0, x_0; d)\|_r \leq \sigma\left(\beta\left(t_0 - \left[\frac{t_0}{T}\right]T\right)\|x_0\|_r, t - t_0\right)$$

Since $0 \leq t_0 - \left[\frac{t_0}{T}\right]T < T$, for all $t_0 \geq 0$, it follows that the following estimate holds for all $(t_0, x_0, d) \in \Re^+ \times C^0([-r,0]; \Re^n) \times M_D$ and $t \in [t_0, +\infty)$:

$$\|\phi(t, t_0, x_0; d)\|_r \leq \tilde{\sigma}\left(\|x_0\|_r, t - t_0\right)$$

where $\tilde{\sigma}(s, t) := \sigma(rs, t)$ and $r := \max\{\beta(t); 0 \leq t \leq T\}$. The previous estimate in conjunction with Definition 2.5 implies that $0 \in C^0([-r,0]; \Re^n)$ is URGAS for system (2.2). The proof is complete. ◁

## 2. III. Differential Inequalities and Dini Derivatives for Functionals

Let $x \in C^0([-r,0]; \Re^n)$. By $E_h(x; v)$, where $0 \leq h < r$ and $v \in \Re^n$ we denote the following operator:

$$E_h(x; v) := \begin{cases} x(0) + (\theta + h)v & \text{for } -h < \theta \leq 0 \\ x(\theta + h) & \text{for } -r \leq \theta \leq -h \end{cases} \tag{2.8}$$

Notice that we have: $\|E_h(x; v) - E_h(x, w)\|_r \leq h|v - w|$, for all $x \in C^0([-r,0]; \Re^n)$, $0 \leq h < r$ and $v, w \in \Re^n$. Let $V : \Re^+ \times C^0([-r,0]; \Re^n) \to \Re$. We define

$$V^0(t,x;v) := \limsup_{\substack{h \to 0^+ \\ y \to 0, y \in C^0([-r,0];\Re^n)}} \frac{V(t+h, E_h(x;v) + hy) - V(t,x)}{h} \qquad (2.9)$$

If there exist constants $L, \delta > 0$ in such a way that $|V(\tau,y) - V(\tau,x)| \leq L\|x-y\|_r$ for all $|\tau - t| + \|y - x\|_r \leq \delta$, then $|V^0(t,x;v) - V^0(t,x;w)| \leq L|v-w|$ for all $v, w \in \Re^n$.

The following lemma presents some elementary properties of this generalized derivative. Notice that the function $(t,x,v) \to V^0(t,x;v)$ may take values in the extended real number system $\Re^* = [-\infty, +\infty]$.

**Lemma 2.7** *Let $V: \Re^+ \times C^0([-r,0]; \Re^n) \to \Re$ and let $x \in C^0([t_0 - r, t_{max}); \Re^n)$ a solution of (2.2) under hypotheses (H1-4) corresponding to certain $d \in M_D$. Then it holds that*

$$\limsup_{h \to 0^+} \frac{V(t+h, T_r(t+h)x) - V(t, T_r(t)x)}{h} \leq V^0(t, T_r(t)x; f(t, T_r(t)x, d(t))), \quad \forall t \in [t_0, t_{max}) \qquad (2.10)$$

**Proof** Inequality (2.10) follows directly from definition (2.9) and the following identity:

$$T_r(t+h)x - E_h(T_r(t)x; f(t, T_r(t)x, d(t))) = \begin{cases} x(t+h+\theta) - x(t) - (\theta + h)\dot{x}(t) & \text{for } -h < \theta \leq 0 \\ 0 & \text{for } -r \leq \theta \leq -h \end{cases} = h y_h$$

where $t \in [t_0, t_{max})$ and

$$y_h := \begin{cases} \frac{\theta + h}{h}\left(\frac{x(t+\theta+h) - x(t)}{\theta + h} - f(t, T_r(t)x, d(t))\right) & \text{for } -h < \theta \leq 0 \\ 0 & \text{for } -r \leq \theta \leq -h \end{cases}$$

with $\|y_h\|_r \leq \sup\left\{\left|\frac{x(t+s) - x(t)}{s} - f(t, T_r(t)x, d(t))\right|; 0 < s \leq h\right\}$. Notice that since $\lim_{h \to 0^+} \frac{x(t+h) - x(t)}{h} = f(t, T_r(t)x, d(t))$ we obtain that $y_h \to 0$ as $h \to 0^+$. The proof is complete. ◁

The following comparison principle is an extension of the comparison principle in [13] and will be used frequently in this paper in conjunction with Lemma 2.7 for the derivation of estimates of values of Lyapunov functionals. Its proof is provided in the Appendix.

**Lemma 2.8 (Comparison Principle)** *Consider the scalar differential equation:*

$$\begin{aligned} \dot{w} &= f(t,w) \\ w(t_0) &= w_0 \end{aligned} \qquad (2.11)$$

*where $f(t,w)$ is continuous in $t \geq 0$ and locally Lipschitz in $w \in J$, where $J \subseteq \Re$ is an open interval. Let $T > t_0$ such that the solution $w(t)$ of the initial value problem (2.11) exists and satisfies $w(t) \in J$ for all $t \in [t_0, T]$. Let $v(t)$ a lower semi-continuous function that satisfies the differential inequality:*

$$D^+ v(t) := \limsup_{h \to 0^+} \frac{v(t+h) - v(t)}{h} \leq f(t, v(t)), \quad \forall t \in [t_0, T) \qquad (2.12)$$

*Suppose furthermore:*

$$v(t_0) \leq w_0 \tag{2.13a}$$

$$v(t) \in J, \forall t \in [t_0, T) \tag{2.13b}$$

*If one of the following holds:*

**(i)** *the mapping $f(t, \cdot)$ is non-decreasing on $J \subseteq \Re$, for each fixed $t \in [t_0, T)$.*

**(ii)** *there exists $\phi \in C^0(\Re^+)$ such that $f(t, w) \leq \phi(t)$, for all $(t, w) \in [t_0, T) \times J$.*

*then $v(t) \leq w(t)$, for all $t \in [t_0, T)$.*

## 2. IV. Statements of Main Results and Examples

We are now in a position to state our main results for non-uniform in time RGAS and URGAS.

**Theorem 2.9** *Consider system (2.2) under hypotheses (H1-4). Then the following statements are equivalent:*

**(a)** $0 \in C^0([-r, 0]; \Re^n)$ *is non-uniformly in time RGAS for (2.2).*

**(b)** *There exists a continuous mapping $(t, x) \in \Re^+ \times C^0([-r, 0]; \Re^n) \to V(t, x) \in \Re^+$, with the following properties:*

**(i)** *There exist functions $a_1, a_2 \in K_\infty$, $\beta \in K^+$ such that:*

$$a_1(\|x\|_r) \leq V(t, x) \leq a_2(\beta(t)\|x\|_r), \forall (t, x) \in \Re^+ \times C^0([-r, 0]; \Re^n) \tag{2.14}$$

**(ii)** *It holds that:*

$$V^0(t, x; f(t, x, d)) \leq -V(t, x), \forall (t, x, d) \in \Re^+ \times C^0([-r, 0]; \Re^n) \times D \tag{2.15}$$
**(infinitesimal decrease property)**

**(iii)** *There exists a non-decreasing function $M : \Re^+ \to \Re^+$ such that for every $R \geq 0$, it holds:*

$$|V(t, y) - V(t, x)| \leq M(R)\|y - x\|_r, \forall t \in [0, R], \forall x, y \in \{x \in C^0([-r, 0]; \Re^n); \|x\|_r \leq R\} \tag{2.16}$$

**(c)** *There exist $\tau \geq 0$, a lower semi-continuous mapping $V : \Re^+ \times C^0([-r-\tau, 0]; \Re^n) \to \Re^+$, constants $R \geq 0$, $c > 0$, functions $a_1, a_2 \in K_\infty$, $\beta_i \in K^+ (i = 1, ..., 4)$ with $\int_0^{+\infty} \beta_4(t)dt = +\infty$, $\mu \in \mathcal{E}$ (see Notations) and $\rho \in C^0(\Re^+; \Re^+)$ being positive definite and locally Lipschitz, such that the following inequalities hold:*

$$a_1(|x(0)|) \leq V(t, x) \leq a_2(\beta_1(t)\|x\|_{r+\tau}), \forall (t, x) \in \Re^+ \times C^0([-r-\tau, 0]; \Re^n) \tag{2.17}$$

$$V^0(t, x; f(t, T_r(0)x, d)) \leq \beta_2(t)V(t, x) + R\beta_3(t), \forall (t, x, d) \in \Re^+ \times C^0([-r-\tau, 0]; \Re^n) \times D \tag{2.18a}$$

$$V^0(t, x; f(t, T_r(0)x, d)) \leq -\beta_4(t)\rho(V(t, x)) + \beta_4(t)\mu\left(\int_0^t \beta_4(s)ds\right), \forall (t, d) \in [\tau, +\infty) \times D, \forall x \in S(t) \tag{2.18b}$$
**(infinitesimal decrease property)**

*where the set-valued map $S(t) \subseteq C^0([-r-\tau, 0]; \Re^n)$ is defined for $t \geq \tau$ by:*

$$S(t) := \left\{ x \in \bar{S}(t); \; x(\theta) = x(-\tau) + \int_{-\tau}^{\theta} f(t+s, T_r(s)x, d(\tau+s))ds, \; \forall \theta \in [-\tau, 0], \; d \in M_D \right\} \quad (2.19a)$$

and $\bar{S}(t) \subseteq C^0([-r-\tau, 0]; \Re^n)$ is any set-valued map satisfying

$$\left\{ x \in C^0([-r-\tau, 0]; \Re^n); \; a_2\left(\beta_1(t)\|x\|_{r+\tau}\right) \geq \eta\left(\int_0^t \beta_4(s)ds, 0, c\right) \right\} \subseteq \bar{S}(t), \; \forall t \geq 0 \quad (2.19b)$$

where $\eta(t, t_0, \eta_0)$ denotes the unique solution of the initial value problem:

$$\dot{\eta} = -\rho(\eta) + \mu(t) \; ; \; \eta(t_0) = \eta_0 \geq 0 \quad (2.19c)$$

**Theorem 2.10** *Consider system (2.2) under hypotheses (H1-4). Then the following statements are equivalent:*

**(a)** $0 \in C^0([-r, 0]; \Re^n)$ *is URGAS for (2.2).*

**(b)** *There exists a continuous mapping* $(t, x) \in \Re^+ \times C^0([-r, 0]; \Re^n) \to V(t, x) \in \Re^+$, *satisfying properties (i), (ii) and (iii) of statement (b) of Theorem 2.9 with* $\beta(t) \equiv 1$. *Moreover, if system (2.2) is* $T$ – *periodic, then* $V$ *is* $T$ – *periodic (i.e.* $V(t+T, x) = V(t, x)$ *for all* $(t, x) \in \Re^+ \times C^0([-r, 0]; \Re^n)$*) and if (2.2) is autonomous then* $V$ *is independent of* $t$.

**(c)** *There exist a lower semi-continuous mapping* $V : \Re^+ \times C^0([-r-\tau, 0]; \Re^n) \to \Re^+$, *functions* $a_1, a_2 \in K_\infty$, $\rho \in C^0(\Re^+; \Re^+)$ *being positive definite and locally Lipschitz and constants* $\beta, \tau \geq 0$ *such that the following inequalities hold:*

$$a_1(|x(0)|) \leq V(t, x) \leq a_2(\|x\|_{r+\tau}), \; \forall (t, x) \in \Re^+ \times C^0([-r-\tau, 0]; \Re^n) \quad (2.20)$$

$$V^0(t, x; f(t, T_r(0)x, d)) \leq \beta V(t, x), \; \forall (t, x, d) \in \Re^+ \times C^0([-r-\tau, 0]; \Re^n) \times D \quad (2.21a)$$

$$V^0(t, x; f(t, T_r(0)x, d)) \leq -\rho(V(t, x)), \; \forall (t, d) \in [\tau, +\infty) \times D, \; \forall x \in S(t) \quad (2.21b)$$
**(infinitesimal decrease property)**

where the set $S(t) \subseteq C^0([-r-\tau, 0]; \Re^n)$ is defined for $t \geq \tau$ by (2.19a) with $\bar{S}(t) := C^0([-r-\tau, 0]; \Re^n)$.

**Remark 2.11:**

**a)** Although the conditions for non-uniform in time RGAS seem more complicated than the corresponding conditions for URGAS, it should be emphasized that the conditions for non-uniform in time RGAS are "weaker" than the corresponding conditions for URGAS. Particularly, the main difference lies in that the infinitesimal decrease condition does not have to be satisfied for states sufficiently close to the equilibrium point in the non-uniform in time case.

**b)** Notice that we demand the infinitesimal decrease property to hold only on a subset of the state space ($S(t)$) which contains the solutions of the system. However, an additional property that guarantees forward completeness on the critical time interval $[t_0, t_0 + \tau]$ has to be satisfied, namely (2.18a) in the non-uniform in time case and (2.21a) in the uniform case. Notice that for finite-dimensional continuous-time systems, it was shown in [1] that this additional property is necessary and sufficient for forward completeness.

**Example 2.12:** Let $b \geq a > 0$, $r \geq 0$ and consider the scalar system:

$$\dot{x}(t) = -d(t)x(t-r)$$
$$x(t) \in \Re, d \in M_D, D := [a,b] \tag{2.22}$$

We will prove (using Theorem 2.10) that $0 \in C^0([-2r,0];\Re)$ is URGAS for (2.22), under the assumption:

$$2b^3r^2 < a \tag{2.23}$$

Clearly, inequality (2.23) is conservative, since it is shown in [8] for the case $d(t) \equiv a$, by using other methods (applicable only to linear systems) that if $2ar < \pi$ then $0 \in C^0([-r,0];\Re)$ is URGAS for (2.22). Here, we consider a class of functionals proposed in [5,21] for systems described by RFDEs and we must make explicit use of the set $S(t)$ involved in statement (c) of Theorem 2.10. Notice that under hypothesis (2.23) there exists $c \in (0,a)$ with $(a-c)(1-2cr) - 2b^3r^2 > 0$. Consider the functional defined on $C^0([-2r,0];\Re)$:

$$V(x) := \frac{1}{2}x^2(0) + \frac{1}{2}\left[(a-c)(1-2cr) - 2b^3r^2\right]\int_{-r}^{0} x^2(s)ds + \frac{1}{2}\left[b^3r + c(a-c)\right]\int_{-2r}^{0}\left(\int_{s}^{0} x^2(l)dl\right)ds \tag{2.24}$$

Since $c \in (0,a)$ with $(a-c)(1-2cr) - 2b^3r^2 > 0$, we conclude that (2.20) is satisfied for this functional with $\tau := r$, $a_1(s) := \frac{1}{2}s^2$ and $a_2(s) := Ks^2$ with $K := \frac{1}{2}(1 + 2r(a-c))$. Moreover, we have:

$$V^0(x;-dx(-r)) = -dx(0)x(-r) + \frac{a-c}{2}x^2(0) - \frac{1}{2}\left[(a-c)(1-2cr) - 2b^3r^2\right]x^2(-r) - \frac{1}{2}\left[b^3r + c(a-c)\right]\int_{-2r}^{0} x^2(l)dl$$

$$\forall (x,d) \in C^0([-2r,0];\Re) \times D \tag{2.25}$$

Completing the squares and using the trivial inequalities $x^2(0) \leq 2V(x)$ and $|d| \leq b$, it follows that (2.21a) is satisfied with $\beta := a - c + \frac{b^2}{(a-c)(1-2cr) - 2b^3r^2}$. Also notice that definition (2.19a) with $\overline{S}(t) := C^0([-2r,0];\Re)$ implies:

$$S(t) := S = \left\{ x \in C^0([-2r,0];\Re); x(\theta) = x(-r) - \int_{-r}^{\theta} \tilde{d}(r+s)x(s-r)ds, \forall \theta \in [-r,0], \tilde{d} \in M_D \right\} \tag{2.26}$$

Combining (2.25) with (2.26) we obtain for all $(x,d) \in S \times D$:

$$V^0(x;-dx(-r)) \leq -dx^2(0) + dx(0)(x(0) - x(-r)) + \frac{a-c}{2}x^2(0) - \frac{1}{2}\left[b^3r + c(a-c)\right]\int_{-2r}^{0} x^2(l)dl$$

$$= \left(\frac{a-c}{2} - d\right)x^2(0) - dx(0)\int_{-r}^{0}\tilde{d}(r+s)x(s-r)ds - \frac{1}{2}\left[b^3r + c(a-c)\right]\int_{-2r}^{0} x^2(l)dl$$

Clearly, we have $\left|dx(0)\int_{-r}^{0}\tilde{d}(r+s)x(s-r)ds\right| \leq \frac{d}{2}x^2(0) + \frac{d}{2}\left|\int_{-r}^{0}\tilde{d}(r+s)x(s-r)ds\right|^2$. Moreover, it holds that $\left|\int_{-r}^{0}\tilde{d}(r+s)x(s-r)ds\right|^2 \leq r\int_{-r}^{0}\tilde{d}^2(r+s)x^2(s-r)ds \leq rb^2\int_{-2r}^{0} x^2(s)ds$ for all $\tilde{d} \in M_D$. The previous inequalities in conjunction with the fact $d \in [a,b]$, give:

$$V^0(x;-dx(-r)) \le -\frac{c}{2}\left(x^2(0)+(a-c)\int_{-2r}^{0} x^2(l)dl\right), \text{ for all } (x,d) \in S \times D$$

On the other hand, since $\int_{-2r}^{0}\left(\int_{s}^{0} x^2(l)dl\right)ds \le 2r\int_{-2r}^{0} x^2(l)dl$ and $\int_{-r}^{0} x^2(s)ds \le \int_{-2r}^{0} x^2(l)dl$, by definition (2.24) we obtain:

$$2V(x) \le x^2(0)+(a-c)\int_{-2r}^{0} x^2(s)ds, \quad \forall x \in C^0([-2r,0];\Re)$$

The two previous inequalities give:

$$V^0(x;-dx(-r)) \le -cV(x), \text{ for all } (x,d) \in S \times D$$

which implies that (2.21b) is also satisfied with $\rho(s) := cs$. We conclude that statement (c) of Theorem 2.10 is satisfied and consequently, $0 \in C^0([-r,0];\Re)$ is URGAS for (2.22). ◁

The following example illustrates the use of statement (c) of Theorem 2.9 for a time-varying nonlinear system described by RFDEs.

**Example 2.13:** Consider the nonlinear planar system:

$$\begin{aligned}
\dot{x}(t) &= -a(t)x(t-1) \\
\dot{y}(t) &= -y(t)+d(t)\exp(t)x^2(t) \\
(x(t),y(t)) &\in \Re^2, \, d \in M_D, \, D := [-1,1], \, t \ge 0
\end{aligned} \quad (2.27)$$

where $a(t) := \begin{cases} 2\sin^2(\pi t) & t \in [2k,2k+1] \\ 0 & t \in (2k-1,2k) \end{cases}$ for each integer $k$. It is shown in [8] (pages 87-88) that the solution of (2.27) satisfies $x(t) = 0$ for all $t \ge t_0 + 4$ and for every initial condition $x_0 \in C^0([-1,0];\Re)$. Here we prove that the equilibrium point $0 \in C^0([-1,0];\Re^2)$ is non-uniformly in time RGAS for (2.27). Consider the Lyapunov functional defined on $C^0([-6,0];\Re^2)$:

$$V(t,x,y) := \frac{1}{2}x^2(0)+\frac{1}{2}\exp(2t)x^4(0)+\int_{-1}^{0}(x^2(s)+x^4(s))ds+\frac{1}{2}y^2(0) \quad (2.28)$$

Notice that inequalities (2.17) hold for the above functional with $r=1$, $\tau = 5$, $a_1(s) = \frac{1}{2}s^2$, $\beta_1(t) := \exp(t)$ and $a_2(s) = 2s^2 + 4s^4$. Defining $f(t,x,y,d) := (-a(t)x(-1),-y(0)+d\exp(t)x^2(0))'$ we obtain:

$$\begin{aligned}
V^0(t,x,y;f(t,x,y,d)) = &-a(t)x(0)x(-1)+(1+\exp(2t))x^4(0)-2a(t)\exp(2t)x^3(0)x(-1) \\
&+x^2(0)-x^2(-1)-x^4(-1)-y^2(0)+d\exp(t)y(0)x^2(0)
\end{aligned} \quad (2.29)$$

Using the Young inequalities $\exp(2t)|x(0)|^3|x(-1)| \le \frac{3}{4}\exp\left(\frac{8t}{3}\right)x^4(0)+\frac{1}{4}x^4(-1)$, $|x(0)||x(-1)| \le \frac{1}{2}x^2(0)+\frac{1}{2}x^2(-1)$ and $\exp(t)|y(0)||x(0)|^2 \le \frac{1}{2}y^2(0)+\frac{1}{2}\exp(2t)x^4(0)$, in conjunction with the fact that $|d| \le 1$ and $|a(t)| \le 2$ for all $t \ge 0$, we obtain:

$$V^0(t,x,y;f(t,x,y,d)) \le 2x^2(0)+\left(1+3\exp\left(\frac{8t}{3}\right)+\frac{3}{2}\exp(2t)\right)x^4(0)-\frac{1}{2}y^2(0), \, \forall(t,x,d) \in \Re^+ \times C^0([-6,0];\Re^2) \times [-1,1]$$

which directly implies (2.18a) with $\beta_2(t) := 12\exp(t)$, $R := 0$ and any $\beta_3 \in K^+$. Let $\bar{S}(t) := C^0([-6,0];\Re^2)$ (and notice that (2.19b) is automatically satisfied for all $c > 0$) and since the solution of (2.27) satisfies $x(t) = 0$ for all $t \geq t_0 + 4$ and for every initial condition $x_0 \in C^0([-1,0];\Re)$, the set-valued map $S(t) \subseteq C^0([-6,0];\Re^2)$ satisfies for all $t \geq 5$:

$$S(t) \subseteq \{(x, y) \in C^0([-6,0];\Re^2) \,;\, x(\theta) = 0, \forall \theta \in [-1,0]\}$$

Consequently, by virtue of definition (2.28) and equality (2.29) we get:

$$V^0(t, x, y; f(t, x, y, d)) \leq -2V(t, x), \quad \forall (t, d) \in [5, +\infty) \times [-1,1], \forall x \in S(t)$$

and thus (2.18b) holds with $\beta_4(t) := 2$, $\rho(s) := s$ and $\mu(t) := 0$. We conclude that statement (c) of Theorem 2.9 is satisfied and consequently, $0 \in C^0([-1,0];\Re^2)$ is non-uniformly in time RGAS for (2.27). ◁

## 3. Proofs of Main Results

The implication (b) $\Rightarrow$ (c) is obvious for both theorems.

**Implication (c) $\Rightarrow$ (a) for Theorem 2.9**

Let arbitrary $(t_0, x_0) \in \Re^+ \times C^0([-r,0];\Re^n)$ and $d \in M_D$. As remarked in Section 2.I there exists $t_{\max} > t_0$ such that the initial-value problem (2.2) with initial condition $T_r(t_0)x = x_0 \in C^0([-r,0];\Re^n)$ has a unique solution $x(t)$ defined on $[t_0 - r, t_{\max})$. Setting $x(t) := x(t_0 - r)$ for $t \in [t_0 - r - \tau, t_0 - r]$, we may assume that for each time $t \in [t_0, t_{\max})$ the unique solution of (2.2) belongs to $C^0([t_0 - r - \tau, t];\Re^n)$. Moreover, we have $\|T_{r+\tau}(t_0)x\|_{r+\tau} = \|T_r(t_0)x\|_r = \|x_0\|_r$.

Let $V(t) := V(t, T_{r+\tau}(t)x)$, which is a lower semi-continuous function on $[t_0, t_{\max})$. Notice that for all $t \in [t_0, t_{\max})$, by virtue of Lemma 2.7 we obtain:

$$D^+V(t) \leq V^0(t, T_{r+\tau}(t)x; f(t, T_r(0)T_{r+\tau}(t)x, d(t))), \text{ for all } t \in [t_0, t_{\max}) \quad (3.1)$$

Inequality (3.1) in conjunction with inequality (2.18a) gives:

$$D^+V(t) \leq \beta_2(t)V(t) + R\beta_3(t), \text{ for all } t \in [t_0, t_{\max}) \quad (3.2)$$

By virtue of Lemma 2.8 (comparison principle, case (i)), we obtain:

$$V(t) \leq \exp\left(\int_{t_0}^t \beta_2(s)ds\right)\left(V(t_0) + R\int_{t_0}^t \beta_3(s)ds\right), \text{ for all } t \in [t_0, t_{\max}) \quad (3.3)$$

The above inequality in conjunction with inequality (2.17) and the fact $\|T_{r+\tau}(t_0)x\|_{r+\tau} = \|T_r(t_0)x\|_r = \|x_0\|_r$, implies that:

$$a_1(|x(t)|) \leq \exp\left(\int_{t_0}^t \beta_2(s)ds\right)\left(a_2(\beta_1(t_0)\|x_0\|_r) + R\int_{t_0}^t \beta_3(s)ds\right), \text{ for all } t \in [t_0, t_{\max}) \quad (3.4)$$

By virtue of the conclusions of Section 2.I, it follows that the solution of (2.2) is defined on $[t_0,+\infty)$ and satisfies (3.3), (3.4) for all $t \in [t_0,+\infty)$. Clearly, inequality (3.4) implies that $0 \in C^0([-r,0];\Re^n)$ is Robustly Forward Complete (see Lemma 2.3). Moreover, by definitions (2.19a,b) it follows that:

$$\text{If } a_2\left(\beta_1(t)\|T_{r+\tau}(t)x\|_{r+\tau}\right) \geq \eta\left(\int_0^t \beta_4(s)ds, 0, c\right) \text{ and } t \geq t_0 + \tau \text{ then } T_{r+\tau}(t)x \in S(t) \qquad (3.5)$$

We proceed by observing the following facts:

**Fact I:** Suppose that $T_{r+\tau}(t)x \in S(t)$ for all $t \in [a,b)$, where $a \geq t_0 + \tau$. Then it holds that

$$D^+V(t) \leq -\beta_4(t)\rho(V(t,T_{r+\tau}(t)x)) + \beta_4(t)\mu\left(\int_0^t \beta_4(s)ds\right), \quad \forall t \in [a,b) \qquad (3.6)$$

This fact can be shown easily using (2.18b), (3.1) and (3.5).

**Fact II:** Suppose that $T_{r+\tau}(t)x \in S(t)$ for all $t \in [a,b)$, where $a \geq t_0 + \tau$. Then the following estimate holds:

$$V(t,T_{r+\tau}(t)x) \leq \eta\left(\int_0^t \beta_4(s)ds, \int_0^a \beta_4(s)ds, V(a,x(a))\right), \quad \forall t \in [a,b) \qquad (3.7)$$

This fact is an immediate consequence of Fact I and Lemma 2.8 (Comparison Principle, case (ii)).

We define the following disjoint sets:

$$A^+ := \left\{ t \in [t_0+\tau,+\infty) ; a_2(\beta_1(t)\|T_{r+\tau}(t)x\|_{r+\tau}) > \eta\left(\int_0^t \beta(s)ds, 0, c\right) \right\} \qquad (3.8)$$

$$A^- := \left\{ t \in [t_0+\tau,+\infty) ; a_2(\beta_1(t)\|T_{r+\tau}(t)x\|_{r+\tau}) \leq \eta\left(\int_0^t \beta(s)ds, 0, c\right) \right\} \qquad (3.9)$$

where $c > 0$ and $\eta(t,t_0,\eta_0)$ are defined in (2.19b,c). Obviously $[t_0+\tau,+\infty) = A^+ \cup A^-$. Notice that by virtue of definitions (2.19b) and (3.5) if $t \in A^+$ then $T_{r+\tau}(t)x \in S(t)$. Moreover, notice that the set $A^+ \setminus \{t_0+\tau\}$ is open. Thus $A^+ \setminus \{t_0+\tau\}$ is either empty or it decomposes into a finite number or a countable infinity of open and disjoint intervals $(a_k,b_k)$ with $a_k < b_k$. When $t_0+\tau \in A^+$ we obviously have the latter case. We distinguish the cases:

**Case A.** $t_0+\tau \notin A^+$ and $A^+ \setminus \{t_0+\tau\}$ is not empty. In this case the set $A^+ \setminus \{t_0+\tau\}$ decomposes into a finite number or a countable infinity of open and disjoint intervals $(a_k,b_k)$ with $a_k < b_k$ for $k=1,\ldots$. Furthermore, by continuity of the solution $T_{r+\tau}(t)x$ it follows that $T_{r+\tau}(a_k)x \in S(a_k)$ and thus $T_{r+\tau}(t)x \in S(t)$ for all $t \in [a_k,b_k)$. Clearly, by Fact II, the following estimate will hold:

$$V(t,T_{r+\tau}(t)x) \leq \eta\left(\int_0^t \beta_4(s)ds, \int_0^{a_k} \beta_4(s)ds, V(a_k,T_{r+\tau}(a_k)x)\right), \quad \forall t \in [a_k,b_k) \qquad (3.10)$$

The fact that $a_k \notin A^+$ implies that $a_k \in A^-$ and consequently by virtue of (2.17) and definition (3.9) we have:

$$V(a_k, T_{r+\tau}(a_k)x) \leq \eta\left(\int_0^{a_k} \beta_4(s)ds, 0, c\right) \quad (3.11)$$

Estimates (3.10) and (3.11) provide the following estimate:

$$V(t, T_{r+\tau}(t)x) \leq \eta\left(\int_0^t \beta_4(s)ds, \int_0^{a_k} \beta_4(s)ds, \eta\left(\int_0^{a_k} \beta_4(s)ds, 0, c\right)\right) = \eta\left(\int_0^t \beta_4(s)ds, 0, c\right), \quad \forall t \in [a_k, b_k) \quad (3.12)$$

When $t \notin [a_k, b_k)$, it follows that $t \in A^-$ and consequently by virtue of definition (3.9) we have:

$$V(t, T_{r+\tau}(t)x) \leq \eta\left(\int_0^t \beta_4(s)ds, 0, c\right), \quad \forall t \notin [a_k, b_k) \quad (3.13)$$

Estimates (3.12) and (3.13) provide the following estimate:

$$V(t, T_{r+\tau}(t)x) \leq \eta\left(\int_0^t \beta_4(s)ds, 0, c\right), \quad \forall t \in [t_0 + \tau, +\infty) \quad (3.14)$$

**Case B.** The set $A^+ \setminus \{t_0 + \tau\}$ is empty. In this case we have $t_0 + \tau \notin A^+$ and consequently it follows that $A^- = [t_0 + \tau, +\infty)$. Therefore by virtue of definition (3.9) we have that estimate (3.14) holds.

**Case C.** $t_0 + \tau \in A^+$ and $A^+ \setminus \{t_0 + \tau\}$ is not empty. In this case there exists a time $b > t_0 + \tau$ and an open set $\widetilde{A}$ such that $A^+ = [t_0 + \tau, b) \cup \widetilde{A}$. For $t \in [t_0 + \tau, b)$ it follows that $T_{r+\tau}(t)x \in S(t)$ and thus by Fact II we obtain the estimate:

$$V(t, T_{r+\tau}(t)x) \leq \eta\left(\int_0^t \beta_4(s)ds, \int_0^{t_0+\tau} \beta_4(s)ds, V(t_0+\tau, T_{r+\tau}(t_0+\tau)x)\right), \quad \forall t \in [t_0+\tau, b) \quad (3.15)$$

For the case $b = +\infty$, the estimate above holds for all $t \in [t_0 + \tau, +\infty)$. For the case $b < +\infty$, we have $b \notin A^+$ and thus we may repeat the analysis in cases A and B for the rest interval.

The analysis above shows that in any case the following estimate holds:

$$V(t, T_{r+\tau}(t)x) \leq \eta\left(\int_0^t \beta_4(s)ds, \int_0^{t_0+\tau} \beta_4(s)ds, V(t_0+\tau, T_{r+\tau}(t_0+\tau)x)\right) + \eta\left(\int_0^t \beta_4(s)ds, 0, c\right), \quad \forall t \in [t_0+\tau, +\infty) \quad (3.16)$$

Lemma 5.2 in [10] implies that there exist a function $\sigma(\cdot) \in KL$ and a constant $M > 0$ such that the following inequalities are satisfied for all $t_0 \geq 0$:

$$0 \leq \eta(t, t_0, \eta_0) \leq \sigma(\eta_0 + M, t - t_0), \quad \forall t \geq t_0, \forall \eta_0 \geq 0 \quad (3.17)$$

Inequalities (3.16) and (3.17) imply that the following estimate holds:

$$V(t) \leq \eta\left(\int_0^t \beta_4(s)ds, \int_0^{t_0+\tau} \beta_4(s)ds, V(t_0+\tau)\right) + \eta\left(\int_0^t \beta_4(s)ds, 0, c\right) \leq 2\sigma\left(V(t_0+\tau) + c + M, \int_{t_0+\tau}^t \beta_4(s)ds\right),$$

$$\forall t \in [t_0 + \tau, +\infty) \quad (3.18)$$

It follows from (2.17), (3.3) and (3.18) that the following estimate holds:

$$a_1(|x(t)|) \leq 2\sigma\left(\exp\left(\int_{t_0}^{t_0+\tau}\beta_2(s)ds\right)\left(a_2(\beta_1(t_0)\|x_0\|_r) + R\int_{t_0}^{t_0+\tau}\beta_3(s)ds\right) + c + M, \int_{t_0+\tau}^{t}\beta_4(s)ds\right), \forall t \in [t_0+\tau,+\infty) \quad (3.19)$$

Estimate (3.19) shows that the Robust Attractivity property is satisfied. Thus, by virtue of Lemma 2.3, the equilibrium point $0 \in C^0([-r,0];\Re^n)$ is non-uniformly in time RGAS for system (2.2).

**Implication (c) $\Rightarrow$ (a) for Theorem 2.10**

Let arbitrary $(t_0, x_0) \in \Re^+ \times C^0([-r,0];\Re^n)$ and $d \in M_D$. As remarked in Section 2.I there exists $t_{\max} > t_0$ such that the initial-value problem (2.2) with initial condition $T_r(t_0)x = x_0 \in C^0([-r,0];\Re^n)$ has a unique solution $x(t)$ defined on $[t_0-r, t_{\max})$. Setting $x(t) := x(t_0-r)$ for $t \in [t_0-r-\tau, t_0-r]$, we may assume that for each time $t \in [t_0, t_{\max})$ the unique solution of (2.2) belongs to $C^0([t_0-r-\tau,t];\Re^n)$. Moreover, we have $\|T_{r+\tau}(t_0)x\|_{r+\tau} = \|T_r(t_0)x\|_r = \|x_0\|_r$.

Let $V(t) := V(t, T_{r+\tau}(t)x)$, which is a lower semi-continuous function on $[t_0, t_{\max})$. Notice that, by virtue of Lemma 2.7, it follows that (3.1) holds for all $t \in [t_0, t_{\max})$. Inequality (3.1) in conjunction with (2.21a) gives:

$$D^+V(t) \leq \beta V(t), \text{ for all } t \in [t_0, t_{\max}) \quad (3.20)$$

By virtue of Lemma 2.8 (comparison principle, case (i)), we obtain:

$$V(t) \leq \exp(\beta(t-t_0))V(t_0), \text{ for all } t \in [t_0, t_{\max}) \quad (3.21)$$

Moreover, by definition (2.19a) with $\overline{S}(t) := C^0([-r-\tau,0];\Re^n)$, it follows that:

$$\text{If } t \geq t_0 + \tau \text{ then } T_{r+\tau}(t)x \in S(t) \quad (3.22)$$

and consequently, by (3.1) and (2.21b), the following differential inequality is satisfied:

$$D^+V(t) \leq -\rho(V(t)), \forall t \in [t_0+\tau,+\infty) \quad (3.23)$$

Let $\eta(t,t_0,\eta_0)$ denotes the unique solution of the initial value problem (2.19c) with $\mu \equiv 0$. Then, by virtue of Lemma 2.8 (comparison principle, case (ii)) the following estimate holds:

$$V(t, T_{r+\tau}(t)x) \leq \eta(t, t_0+\tau, V(t_0+\tau)), \forall t \in [t_0+\tau,+\infty) \quad (3.24)$$

Lemma 4.4 in [19] implies that there exists $\sigma(\cdot) \in KL$ such that the following inequalities are satisfied for all $t_0 \geq 0$:

$$0 \leq \eta(t,t_0,\eta_0) \leq \sigma(\eta_0, t-t_0), \forall t \geq t_0, \forall \eta_0 \geq 0 \quad (3.25)$$

Combining (3.21), (3.24) and (3.25) we obtain:

$$V(t, T_{r+\tau}(t)x) \leq \tilde{\sigma}(V(t_0), t-t_0), \forall t \in [t_0,+\infty) \quad (3.26)$$

where $\tilde{\sigma}(s,t) := \sigma(\exp(\beta\tau)s, t-\tau)$ for $t \geq \tau$ and $\tilde{\sigma}(s,t) := \exp(\tau-t)\sigma(\exp(\beta\tau)s, 0)$ for $0 \leq t < \tau$. Let the $KL$ function $\overline{\sigma}(s,t) := a_1^{-1}(\sigma(a_2(s), t))$ and notice that inequalities (2.20) and (3.16) imply the following estimate:

$$|x(t)| \leq \overline{\sigma}(\|x_0\|_r, t-t_0), \forall t \geq t_0 \quad (3.27)$$

It follows from (3.27) that (2.7) holds for the $KL$ function $\sigma(s,t) := \bar{\sigma}(s, t-r)$ for $t \geq r$ and $\sigma(s,t) := \exp(r-t)\bar{\sigma}(s, 0)$ for $0 \leq t < r$.

**Implication (a) $\Rightarrow$ (b) for Theorem 2.9**

The analysis followed here is similar to the corresponding analysis in [3] for finite-dimensional continuous-time systems.

Since $0 \in C^0([-r,0];\mathfrak{R}^n)$ is non-uniformly in time RGAS for (2.2), there exist functions $\sigma \in KL$, $\beta \in K^+$ such that estimate (2.6) holds for all $(t_0, x_0, d) \in \mathfrak{R}^+ \times C^0([-r,0];\mathfrak{R}^n) \times M_D$ and $t \in [t_0, +\infty)$. Moreover, by recalling Proposition 7 in [28] there exist functions $\tilde{a}_1$, $\tilde{a}_2$ of class $K_\infty$, such that the $KL$ function $\sigma(s,t)$ is dominated by $\tilde{a}_1^{-1}(\exp(-2t)\tilde{a}_2(s))$. Thus, by taking into account estimate (2.6), we have:

$$\tilde{a}_1\left(\|\phi(t,t_0,x_0;d)\|_r\right) \leq \exp(-2(t-t_0))\tilde{a}_2\left(\beta(t_0)\|x_0\|_r\right), \forall t \geq t_0 \geq 0, x_0 \in C^0([-r,0];\mathfrak{R}^n), d \in M_D \tag{3.28}$$

Without loss of generality we may assume that $\tilde{a}_1 \in K_\infty$ is globally Lipschitz on $\mathfrak{R}^+$ with unit Lipschitz constant, namely, $|\tilde{a}_1(s_1) - \tilde{a}_1(s_2)| \leq |s_1 - s_2|$ for all $s_1, s_2 \geq 0$. To see this notice that we can always replace $\tilde{a}_1 \in K_\infty$ by the function $\bar{a}_1(s) := \inf\left\{\min\left\{\frac{1}{2}y, a(y)\right\} + |y-s|; y \geq 0\right\}$, which is of class $K_\infty$, globally Lipschitz on $\mathfrak{R}^+$ with unit Lipschitz constant and satisfies $\bar{a}_1(s) \leq \tilde{a}_1(s)$. Moreover, without loss of generality we may assume that $\beta \in K^+$ is non-decreasing.

Making use of (2.5) and (3.28), we obtain the following elementary property for the solution of (2.2):

$$\|\phi(t,t_0,x;d) - \phi(t,t_0,y;d)\|_r \leq \exp\left(\tilde{L}(t,\|x\|_r + \|y\|_r)(t-t_0)\right)\|x-y\|_r \tag{3.29}$$

for all $t \geq t_0$ and $(t_0, x, y, d) \in \mathfrak{R}^+ \times C^0([-r,0];\mathfrak{R}^n) \times C^0([-r,0];\mathfrak{R}^n) \times M_D$

where

$$\tilde{L}(t,s) := L(t, 2\tilde{a}_1^{-1}(\tilde{a}_2(\beta(t)s)))$$

and $L(\cdot)$ is the function involved in (2.3). Furthermore, under hypotheses (H1-4), Lemma 3.2 in [12] implies the existence of functions $\zeta \in K_\infty$ and $\gamma \in K^+$ such that:

$$|f(t,x,d)| \leq \zeta\left(\gamma(t)\|x\|_r\right), \forall (t,x,d) \in \mathfrak{R}^+ \times C^0([-r,0];\mathfrak{R}^n) \times D$$

Without loss of generality, we may assume that $\gamma \in K^+$ is non-decreasing. Since $x(t) = x(0) + \int_{t_0}^{t} f(\tau, T_r(\tau)x, d(\tau))d\tau$, using the previous inequality in conjunction with (3.28) we obtain:

$$|\phi(0,t,t_0,x;d) - x(0)| \leq (t-t_0)G_1(t,\|x\|_r)$$
$$G_1(t,s) := \zeta\left(\gamma(t)\tilde{a}_1^{-1}(\tilde{a}_2(\beta(t)s))\right)$$

and consequently

$$\|\phi(t,t_0,x;d) - x\|_r \leq (t-t_0)G_1(t,\|x\|_r) + G_2(x, t-t_0) \tag{3.30}$$

for all $t \geq t_0$ and $(t_0, x, d) \in \mathfrak{R}^+ \times C^0([-r,0];\mathfrak{R}^n) \times M_D$

where the functional

$$G_2(x,h) := \sup\{|x(0)-x(\theta)|; \theta \in [-\min(h,r),0]\} + \begin{cases} 0 & if \quad h \geq r \\ \sup\{|x(\theta+h)-x(\theta)|; \theta \in [-r,-h]\} & if \quad 0 \leq h < r \end{cases}$$

is defined for all $(x,h) \in C^0([-r,0]; \Re^n) \times \Re^+$. Notice that $\lim_{h \to 0^+} G_2(x,h) = 0$ for all $x \in C^0([-r,0]; \Re^n)$ and consequently for every $\varepsilon > 0$, $R \geq 0$, $x \in C^0([-r,0]; \Re^n)$, there exists $T(\varepsilon, R, x) > 0$ such that:

$$t_0 \leq t \leq t_0 + T(\varepsilon, R, x) \Rightarrow \|\phi(t, t_0, x; d) - x\|_r \leq \varepsilon \quad (3.31)$$
$$\text{for all } (t_0, x, d) \in [0, R] \times C^0([-r,0]; \Re^n) \times M_D$$

We define for all $q \in Z^+$:

$$U_q(t,x) := \sup\{\max\{0, \tilde{a}_1(\|\phi(\tau,t,x;d)\|_r) - q^{-1}\}\exp((\tau-t)) : \tau \geq t, d \in M_D\} \quad (3.32)$$

Clearly, estimate (2.6) and definition (3.32) imply that:

$$\max\{0, \tilde{a}_1(\|x\|_r) - q^{-1}\} \leq U_q(t,x) \leq \tilde{a}_2(\beta(t)\|x\|_r), \quad \forall (t,x,q) \in \Re^+ \times C^0([-r,0]; \Re^n) \times Z^+ \quad (3.33)$$

Moreover, by definition (3.32) we obtain for all $(h,t,x,d,q) \in \Re^+ \times \Re^+ \times C^0([-r,0]; \Re^n) \times M_D \times Z^+$:

$$U_q(t+h, \phi(t+h,t,x;d)) \leq \exp(-h)U_q(t,x) \quad (3.34)$$

By virtue of estimate (2.6) it follows that for every $(q,R) \in Z^+ \times \Re^+$, $\tau \geq t + \tilde{T}(R,q)$, $(t,d) \in [0,R] \times M_D$, and $x \in C^0([-r,0]; \Re^n)$ with $\|x\|_r \leq R$, it holds: $\tilde{a}_1(\|\phi(\tau,t,x;d)\|_r) \leq \exp(-2(\tau-t))\tilde{a}_2(\beta(t)\|x\|_r) \leq q^{-1}$, where

$$\tilde{T}(R,q) := \max\left\{0, \frac{1}{2}\log(q\tilde{a}_2(\beta(R)R))\right\} \quad (3.35)$$

Thus, by virtue of definition (3.32), we conclude that:

$$U_q(t,x) = \sup\{\max\{0, \tilde{a}_1(\|\phi(\tau,t,x;d)\|_r) - q^{-1}\}\exp((\tau-t)) : t \leq \tau \leq t+\xi, d \in M_D\}$$
$$\text{for all } \xi \geq \tilde{T}(\max\{t, \|x\|_r\}, q), x \neq 0 \quad (3.36)$$

It follows by taking into account (3.36) that for all $t \in [0, R]$, and $(x,y) \in C^0([-r,0]; \Re^n) \times C^0([-r,0]; \Re^n)$ with $\|x\|_r \leq R$, $\|y\|_r \leq R$, it holds:

$$|U_q(t,y) - U_q(t,x)| =$$
$$|\sup\{\max\{0, \tilde{a}_1(\|\phi(\tau,t,y;d)\|_r) - q^{-1}\}\exp((\tau-t)) : t \leq \tau \leq t+\tilde{T}(R,q), d \in M_D\}$$
$$- \sup\{\max\{0, \tilde{a}_1(\|\phi(\tau,t,x;d)\|_r) - q^{-1}\}\exp((\tau-t)) : t \leq \tau \leq t+\tilde{T}(R,q), d \in M_D\}|$$
$$\leq \sup\{\exp((\tau-t))|\tilde{a}_1(\|\phi(\tau,t,y;d)\|_r) - \tilde{a}_1(\|\phi(\tau,t,y;d)\|_r)| : t \leq \tau \leq t+\tilde{T}(R,q), d \in M_D\}$$
$$\leq \sup\{\exp((\tau-t))\|\phi(\tau,t,y;d) - \phi(\tau,t,x;d)\|_r : t \leq \tau \leq t+\tilde{T}(R,q), d \in M_D\}$$
$$(3.37)$$

Notice that in the above inequalities we have used the facts that the functions $\max\{0, s-q^{-1}\}$ and $\tilde{a}_1(s)$ are globally Lipschitz on $\Re^+$ with unit Lipschitz constant. From (3.29) and (3.37) we deduce for all $t \in [0, R]$, and $(x,y) \in C^0([-r,0]; \Re^n) \times C^0([-r,0]; \Re^n)$ with $\|x\|_r \leq R$, $\|y\|_r \leq R$:

$$|U_q(t,y) - U_q(t,x)| \leq G_3(R,q) \|y-x\|_r \tag{3.38}$$

where

$$G_3(R,q) := \exp\left(\widetilde{T}(R,q)\left(1 + \widetilde{L}(R + \widetilde{T}(R,q), 2R)\right)\right) \tag{3.39}$$

Next, we establish continuity with respect to $t$ on $\Re^+ \times C^0([-r,0]; \Re^n)$. Let $R \geq 0$, $q \in Z^+$ arbitrary, $t_1, t_2 \in [0, R]$ with $t_1 \leq t_2$, and $x \in C^0([-r,0]; \Re^n)$ with $\|x\|_r \leq R$. Clearly, we have for all $d \in M_D$:

$$\begin{aligned}|U_q(t_1,x) - U_q(t_2,x)| &\leq (1-\exp(-(t_2-t_1)))U_q(t_1,x) + |\exp(-(t_2-t_1))U_q(t_1,x) - U_q(t_2, \phi(t_2,t_1,x;d))| \\ &+ |U_q(t_2, \phi(t_2,t_1,x;d)) - U_q(t_2,x)|\end{aligned}$$

By virtue of (3.31), (3.34), (3.38) and the previous inequality we obtain for all $t_1, t_2 \in [0, R]$ with $t_1 \leq t_2 \leq t_1 + T(1, R, x)$ (where $T(\varepsilon, R, x) > 0$ is involved in (3.31)) and $d \in M_D$:

$$\begin{aligned}|U_q(t_1,x) - U_q(t_2,x)| &\leq (t_2-t_1)U_q(t_1,x) + \exp(-(t_2-t_1))U_q(t_1,x) - U_q(t_2, \phi(t_2,t_1,x;d)) \\ &+ G_3(R+1,q)\left[G_2(x, t_2-t_1) + (t_2-t_1)G_1(R,R)\right]\end{aligned} \tag{3.40}$$

Definition (3.32) implies that for every $\varepsilon > 0$, there exists $d_\varepsilon \in M_D$ with the following property:

$$U_q(t_1,x) - \varepsilon \leq \sup\left\{\max\left\{0, \widetilde{a}_1\left(\|\phi(\tau,t_1,x;d_\varepsilon)\|_r\right) - q^{-1}\right\}\exp((\tau-t_1)); \tau \geq t_1\right\} \leq U_q(t_1,x) \tag{3.41}$$

Thus using definition (3.32) we obtain:

$$\begin{aligned}&\exp(-(t_2-t_1))U_q(t_1,x) - U_q(t_2, \phi(t_2,t_1,x;d_\varepsilon)) \\ &\leq \max\{A_q(t_1,t_2,x), B_q(t_1,t_2,x)\} - B_q(t_1,t_2,x) + \varepsilon \exp(-(t_2-t_1))\end{aligned} \tag{3.42a}$$

where

$$\begin{aligned}A_q(t_1,t_2,x) &:= \sup\left\{\max\left\{0, \widetilde{a}_1\left(\|\phi(\tau,t_1,x;d_\varepsilon)\|_r\right) - q^{-1}\right\}\exp((\tau-t_2)); t_2 \geq \tau \geq t_1\right\} \\ B_q(t_1,t_2,x) &:= \sup\left\{\max\left\{0, \widetilde{a}_1\left(\|\phi(\tau,t_1,x;d_\varepsilon)\|_r\right) - q^{-1}\right\}\exp((\tau-t_2)); \tau \geq t_2\right\}\end{aligned} \tag{3.42b}$$

Since the functions $\max\{0, s - q^{-1}\}$ and $\widetilde{a}_1(s)$ are globally Lipschitz on $\Re^+$ with unit Lipschitz constant, we obtain:

$$\begin{aligned}&A_q(t_1,t_2,x) - B_q(t_1,t_2,x) \\ &\leq \sup\left\{\max\left\{0, \widetilde{a}_1\left(\|\phi(\tau,t_1,x;d_\varepsilon)\|_r\right) - q^{-1}\right\}\exp((\tau-t_2)); t_2 \geq \tau \geq t_1\right\} - \max\left\{0, \widetilde{a}_1\left(\|\phi(t_2,t_1,x;d_\varepsilon)\|_r\right) - q^{-1}\right\} \\ &\leq \sup\left\{\max\left\{0, \widetilde{a}_1\left(\|\phi(\tau,t_1,x;d_\varepsilon)\|_r\right) - q^{-1}\right\}; t_2 \geq \tau \geq t_1\right\} - \max\left\{0, \widetilde{a}_1\left(\|\phi(t_2,t_1,x;d_\varepsilon)\|_r\right) - q^{-1}\right\} \\ &\leq \sup\left\{\left|\widetilde{a}_1\left(\|\phi(\tau,t_1,x;d_\varepsilon)\|_r\right) - \widetilde{a}_1\left(\|\phi(t_2,t_1,x;d_\varepsilon)\|_r\right)\right|; t_2 \geq \tau \geq t_1\right\} \\ &\leq \sup\left\{\|\phi(\tau,t_1,x;d_\varepsilon) - \phi(t_2,t_1,x;d_\varepsilon)\|_r; t_2 \geq \tau \geq t_1\right\}\end{aligned} \tag{3.43}$$

Notice that by virtue of (3.30), we obtain for all $\tau \in [t_1, t_2]$:

$$\begin{aligned}&\|\phi(\tau,t_1,x;d_\varepsilon) - \phi(t_2,t_1,x;d_\varepsilon)\|_r \\ &\leq \|\phi(\tau,t_1,x;d_\varepsilon) - x\|_r + \|\phi(t_2,t_1,x;d_\varepsilon) - x\|_r \\ &\leq 2(t_2-t_1)G_1(R,R) + 2\sup\{G_2(x,h); h \in [0, t_2-t_1]\}\end{aligned} \tag{3.44}$$

Distinguishing the cases $A_q(t_1,t_2,x) \geq B_q(t_1,t_2,x)$ and $A_q(t_1,t_2,x) \leq B_q(t_1,t_2,x)$ it follows from (3.42a,b), (3.43) and (3.44) that:

$$\exp(-(t_2-t_1))U_q(t_1,x) - U_q(t_2, \phi(t_2,t_1,x;d_\varepsilon)) \leq 2(t_2-t_1)G_1(R,R) + 2\sup\{G_2(x,h); h \in [0, t_2-t_1]\} + \varepsilon$$

Combining the previous inequality with (3.40) we obtain:

$$\begin{aligned}&|U_q(t_1,x)-U_q(t_2,x)|\\&\leq (t_2-t_1)U_q(t_1,x)+\left(2+G_3(R+1,q)\right)\left[\sup\{G_2(x,h)\,;\,h\in[0,t_2-t_1]\}+(t_2-t_1)G_1(R,R)\right]+\varepsilon\end{aligned} \quad (3.45)$$

Since (3.45) holds for all $\varepsilon > 0$, $R \geq 0$, $q \in Z^+$, $x \in C^0([-r,0];\Re^n)$ with $\|x\|_r \leq R$ and $t_1, t_2 \in [0, R]$ with $t_1 \leq t_2 \leq t_1 + T(1, R, x)$, it follows that:

$$|U_q(t_1,x)-U_q(t_2,x)|\leq |t_2-t_1|U_q(t_1,x)+\left(2+G_3(R+1,q)\right)\left[\sup\{G_2(x,h)\,;\,h\in[0,|t_2-t_1|]\}+|t_2-t_1|G_1(R,R)\right]$$

for all $R \geq 0$, $q \in Z^+$, $x \in C^0([-r,0];\Re^n)$ with $\|x\|_r \leq R$ and $t_1, t_2 \in [0, R]$ with $|t_2 - t_1| \leq T(1, R, x)$ (3.46)

Finally, we define:

$$V(t,x):=\sum_{q=1}^{\infty}\frac{2^{-q}U_q(t,x)}{1+G_3(q,q)+(2+G_3(q+1,q))(1+G_1(q,q))} \quad (3.47)$$

Inequality (3.33) in conjunction with definition (3.47) implies (2.14) with $a_2 = \tilde{a}_2$ and

$$a_1(s):=\sum_{q=1}^{\infty}\frac{2^{-q}\max\{0,\,a_1(s)-q^{-1}\}}{1+G_3(q,q)+(2+G_3(q+1,q))(1+G_1(q,q))},$$

which is a function of class $K_\infty$. Moreover, by definition (3.47) and inequality (3.34) we obtain for all $(h,t,x,d) \in \Re^+ \times \Re^+ \times C^0([-r,0];\Re^n) \times M_D$:

$$V(t+h,\phi(t+h,t,x;d))\leq \exp(-h)V(t,x) \quad (3.48)$$

Next define

$$M(R):=1+\sum_{q=1}^{[R]}\frac{2^{-q}G_3(R,q)}{1+G_3(q,q)} \quad (3.49)$$

which is a positive non-decreasing function. Using (3.38) and definition (3.47) as well as the fact $G_3(R,q) \leq G_3(q,q)$ for $q > R$, we may establish (2.16). Finally, by virtue of (3.45), (3.47) and the facts $G_3(R,q) \leq G_3(q,q)$, $G_1(R,R) \leq G_1(q,q)$ for $q > R$, we obtain:

$$|V(t_1,x)-V(t_2,x)|\leq |t_2-t_1|P_1(R)+P_2(R)\sup\{G_2(x,h)\,;\,h\in[0,|t_2-t_1|]\}$$

for all $R \geq 0$, $x \in C^0([-r,0];\Re^n)$ with $\|x\|_r \leq R$ and $t_1, t_2 \in [0, R]$ with $|t_2 - t_1| \leq T(1, R, x)$ (3.50)

where

$$P_1(R):=a_2(\beta(R)R)+1+\sum_{q=1}^{[R]}\frac{2^{-q}(2+G_3(R+1,q))G_1(R,R)}{1+G_3(q,q)+(2+G_3(q+1,q))(1+G_1(q,q))}$$

$$P_2(R):=1+\sum_{q=1}^{[R]}\frac{2^{-q}(2+G_3(R+1,q))}{1+G_3(q,q)+(2+G_3(q+1,q))(1+G_1(q,q))}$$

are positive non-decreasing functions. Inequality (3.50) in conjunction with the fact that $\lim_{h\to 0^+} G_2(x,h) = 0$ for all $x \in C^0([-r,0];\Re^n)$, establishes continuity of $V$ with respect to $t$ on $\Re^+ \times C^0([-r,0];\Re^n)$. Let $d \in D$ and define $\tilde{d}(t) \equiv d$. By definition (2.9) and inequality (3.48), we have for all $(t,x) \in \Re^+ \times C^0([-r,0];\Re^n)$:

$$V^0(t,x;f(t,x,d)) := \limsup_{\substack{h \to 0^+ \\ y \to 0, y \in C^0([-r,0];\Re^n)}} \frac{V(t+h, E_h(x;f(t,x,d))+hy) - V(t,x)}{h}$$

$$\leq \limsup_{h \to 0^+} \frac{V(t+h, \phi(t+h,t,x;\tilde{d})) - V(t,x)}{h} + \limsup_{\substack{h \to 0^+ \\ y \to 0, y \in C^0([-r,0];\Re^n)}} \frac{V(t+h, E_h(x;f(t,x,d))+hy) - V(t+h, \phi(t+h,t,x;\tilde{d}))}{h}$$

$$\leq -V(t,x) + \limsup_{\substack{h \to 0^+ \\ y \to 0, y \in C^0([-r,0];\Re^n)}} \frac{V(t+h, E_h(x;f(t,x,d))+hy) - V(t+h, \phi(t+h,t,x;\tilde{d}))}{h}$$

Let $R \geq \max\{t, \|x\|_r\}$. By definition (2.8) and property (3.31) it follows that $t+h \leq R+1$, $\|\phi(t+h,t,x;\tilde{d})\|_r \leq R+1$, $\|E_h(x;f(t,x,d))+hy\|_r \leq R+1$ for $h$ and $\|y\|_r$ sufficiently small. Using (2.16) we obtain:

$$V^0(t,x;f(t,x,d)) \leq -V(t,x) + M(R+1) \limsup_{h \to 0^+} \frac{\|E_h(x;f(t,x,d)) - \phi(t+h,t,x;\tilde{d})\|_r}{h}$$

As in the proof of Lemma 2.7 we may establish that $\limsup_{h \to 0^+} \frac{\|E_h(x;f(t,x,d)) - \phi(t+h,t,x;\tilde{d})\|_r}{h} = 0$ and consequently the previous inequality shows that (2.15) also holds.

**Implication (a) $\Rightarrow$ (b) for Theorem 2.10**

Since $0 \in C^0([-r,0];\Re^n)$ is URGAS for (2.2), there exist a function $\sigma \in KL$ such that estimate (2.6) holds for all $(t_0, x_0, d) \in \Re^+ \times C^0([-r,0];\Re^n) \times M_D$ and $t \in [t_0, +\infty)$ with $\beta(t) \equiv 1$. Thus all the previous arguments may be repeated for the special case of the constant function $\beta(t) \equiv 1$. We finish the proof with some remarks for the following particular cases:

∗ If (2.2) is $T-periodic$ then for all $(t_0, x_0, d) \in \Re^+ \times C^0([-r,0];\Re^n) \times M_D$ it holds that $\phi(t,t_0,x_0;d) = \phi\left(t - \left[\frac{t_0}{T}\right]T, t_0 - \left[\frac{t_0}{T}\right]T, x_0; P(t_0)d\right)$, where $\left[\frac{t_0}{T}\right]$ denotes the integer part of $\frac{t_0}{T}$ and $P(t_0)d \in M_D$ is defined by:

$$(P(t_0)d)(t) := d\left(t + \left[\frac{t_0}{T}\right]T\right), \quad \forall t \geq 0$$

It follows from definition (3.32) for all $q \in Z^+$:

$$U_q(t,x) := \sup\left\{ \max\left\{0, \tilde{a}_1\left(\left\|\phi\left(\tau - \left[\frac{t}{T}\right]T, t - \left[\frac{t}{T}\right]T, x; P(t)d\right)\right\|_r\right) - q^{-1}\right\} \exp((\tau - t)) : \tau \geq t, d \in M_D \right\}$$

Moreover, since $P(t)M_D = M_D$ for all $t \geq 0$, it follows that $U_q(t,x) = U_q\left(t - \left[\frac{t}{T}\right]T, x\right)$, which directly implies that $U_q$ is $T-periodic$ for all $q \in Z^+$. Definition (3.47) implies that $V$ is $T-periodic$.

∗ If (2.2) is autonomous then (2.2) is $T-periodic$ for all $T > 0$. Consequently, $V$ is $T-periodic$ for all $T > 0$ and thus $V(T,x) = V(0,x)$ for all $T > 0$. Thus $V(t,x) = V(0,x)$ for all $t \geq 0$ and consequently $V$ is independent of $t$.

The proof is complete. ◁

## 5. Conclusions

We establish Lyapunov characterizations for the concepts of *non-uniform in time and uniform robust global asymptotic stability* (RGAS) for uncertain systems described by RFDEs, completely analogous to the corresponding characterizations for continuous-time finite-dimensional uncertain systems, which overcome the limitations imposed by previous works. Particularly, our Lyapunov characterizations apply

- to systems with disturbances that take values in a (not necessarily compact) given set
- to systems described by RFDEs with right-hand sides which are not necessarily bounded with respect to time

The established Lyapunov-like conditions demand the infinitesimal decrease property to hold only on subsets of the state space along with an additional property that guarantees forward completeness. Illustrating examples are also provided.

**Acknowledgements:** The author would like to thank P.A. Bliman and F. Mazenc for our fruitful discussions on systems described by RFDEs.

**References**
[1] Angeli, D. and E.D. Sontag, "Forward Completeness, Unbounded Observability and their Lyapunov Characterizations", *Systems and Control Letters*, 38(4-5), 1999, 209-217.
[2] Aubin, J.P., "Viability Theory", Birkhauser, Boston, 1991.
[3] Bacciotti, A. and L. Rosier, "Liapunov Functions and Stability in Control Theory", Lecture Notes in Control and Information Sciences, 267, Springer-Verlag, London, 2001.
[4] Bacciotti, A., "Stabilization by Means of State Space Switching Rules", Systems and Control Letters, 53(3-4), 2004, 195-201.
[5] Bliman, P.A., "Lyapunov Equation for the Stability of Linear Delay Systems of Retarded and Neutral Type", *IEEE Transactions on Automatic Control*, 47(2), 2002, 327-335.
[6] Burton, T.A. and L. Hatvani, "Stability Theorems for Nonautonomous Functional Differential Equations by Liapunov Functionals", *Tohoku Math. J.*, 41, 1989, 65-104.
[7] Burton, T.A., "Liapunov Functionals, Fixed Points and Stability by Krasnoselskii's Theorem", *Nonlinear Studies*, 9, 2001, 181-190.
[8] Hale, J.K. and S.M.V. Lunel, "Introduction to Functional Differential Equations", Springer-Verlag, New York, 1993.
[9] Jankovic, M., "Control Lyapunov-Razumikhin Functions and Robust Stabilization of Time Delay Systems", *IEEE Transactions on Automatic Control*, 46(7), 2001, 1048-1060.
[10] Karafyllis, I. and J. Tsinias, "A Converse Lyapunov Theorem for Non-Uniform in Time Global Asymptotic Stability and its Application to Feedback Stabilization", *SIAM Journal on Control and Optimization*, 42(3), 2003, 936-965.
[11] Karafyllis, I. and J. Tsinias, "Non-Uniform in Time ISS and the Small-Gain Theorem", *IEEE Transactions on Automatic Control*, 49(2), 2004, 196-216.
[12] Karafyllis, I., "The Non-Uniform in Time Small-Gain Theorem for a Wide Class of Control Systems with Outputs", *European Journal of Control*, 10(4), 2004, 307-323.
[13] Khalil, H.K., "Nonlinear Systems", 2$^{nd}$ Edition, Prentice-Hall, 1996.
[14] Kharitonov, V.L. and D. Melchor-Aguilar, "On Delay Dependent Stability Conditions for Time-Varying Systems", *Systems and Control Letters*, 46, 2002, 173-180.
[15] Kharitonov, V.L., "Lyapunov-Krasovskii Functionals for Scalar Time Delay Equations", *Systems and Control Letters*, 51, 2004, 133-149.
[16] Krasovskii, N.N., "Stability of Motion", Stanford: Stanford University Press, 1963.
[17] Lakshmikantham, V., S. Leela and A.A. Martynyuk, "Stability Analysis of Nonlinear Systems", Marcel Dekker Inc., 1989.
[18] Lakshmikantham, V. and M.R.M. Rao, "Theory of Integrodifferential Equations", Gordon and Breach Science Publishers, U.K., 1995.
[19] Lin, Y., E.D. Sontag and Y. Wang, "A Smooth Converse Lyapunov Theorem for Robust Stability", *SIAM Journal on Control and Optimization*, 34, 1996, 124-160.
[20] Mazenc, F. S. Mondie and S.I. Niculescu, "Global Asymptotic Stabilization for Chains of Integrators with a Delay in the Input", *IEEE Transactions on Automatic Control*, 48(1), 2003, 57-63.
[21] Mazenc, F. and P.A. Bliman, "Backstepping Design for Time-Delay Nonlinear Systems", *Proceedings of the 42$^{nd}$ IEEE Conference on Decision and Control,* Maui, Hawaii, U.S.A., December 2003.
[22] Mazenc, F. and C. Prieur, "Switching Using Distributed Delays", *Proceedings of the 4$^{th}$ IFAC Workshop on Time Delay Systems*, Rocquencourt, France, 2003.

[23] Melchor-Aguilar, D., V. Kharitonov and R. Lozano, "Lyapunov-Krasovskii Functionals for Integral Delay Equations", *Proceedings of the 4th IFAC Workshop on Time Delay Systems*, Rocquencourt, France, 2003.
[24] Morin, P. and C. Samson, "Robust Point-Stabilization of Nonlinear Affine Control Systems", in "Stability and Stabilization of Nonlinear Systems", D. Aeyels, F. Lamnabhi-Lagarrigue and A. van der Schaft (Eds), Springer-Verlag, London, 1999, 215-237.
[25] Niculescu, S.I., "Delay Effects on Stability, A Robust Control Approach", Heidelberg, Germany, Springer-Verlag, 2001.
[26] Savkin, A.V. and R.J. Evans, "Hybrid Dynamical Systems", Birkhauser, Boston, 2002.
[27] Sontag, E.D., "Smooth Stabilization Implies Coprime Factorization", *IEEE Trans. Automat. Contr.*, 34, 1989, 435-443.
[28] Sontag, E.D., "Comments on Integral Variants of ISS", *Systems and Control Letters*, 34, 1998, 93-100.
[29] Teel, A.R., "Connections between Razumikhin-Type Theorems and the ISS Nonlinear Small Gain Theorem", *IEEE Transactions on Automatic Control*, 43(7), 1998, 960-964.
[30] Vorotnikov, V.I., "Partial Stability and Control", Birkhauser, Boston, 1998.
[31] Yoshizawa, T., "Asymptotic Behavior of Solutions in Nonautonomous Systems", in "Trends in Theory and Practice of Nonlinear Differential Equations", V. Lakshmikantham (Eds), Dekker, New York, 1984, 553-562.

# Appendix

**Proof of Lemma 2.8 (Comparison Principle):** It suffices to show that $v(t) \leq w(t)$ on any compact interval $[t_0, t_1] \subset [t_0, T)$. Consider the scalar differential equation:

$$\dot{z} = f(t, z) + \lambda$$
$$z(t_0) = w_0 \tag{A1}$$

where $\lambda$ is a positive constant. On any compact interval $[t_0, t_1] \subset [t_0, T)$, we conclude from Theorem 2.6 in [13] that for every $\varepsilon > 0$ there exists $\delta > 0$ such that if $0 < \lambda < \delta$ then (A1) has a unique solution $z(t, \lambda)$ defined on $[t_0, t_1]$ and satisfies:

$$z(t, \lambda) \in J, \ |z(t, \lambda) - w(t)| < \varepsilon, \ \forall t \in [t_0, t_1] \tag{A2}$$

**Fact I:** $v(t) \leq z(t, \lambda)$, for all $t \in [t_0, t_1)$.

This fact is shown by contradiction. Suppose that there exists $t \in (t_0, t_1)$ such that $v(t) - z(t, \lambda) > 0$. Let the lower semi-continuous function $m(t) := v(t) - z(t, \lambda)$ and define the set:

$$A^+ := \{\tau \in (t_0, t_1) : m(\tau) > 0\} \tag{A3}$$

which by assumption is non-empty. Lower semi-continuity of $m(t) := v(t) - z(t, \lambda)$ implies that $A^+$ is open. Let

$$\tilde{t} := \inf\{t \in A^+\} \tag{A4}$$

Since $A^+$ is open we conclude that $\tilde{t} \notin A^+$, or equivalently that $v(\tilde{t}) \leq z(\tilde{t}, \lambda)$. On the other hand by definition (A4) there exists a sequence $\{\tau_i \in A^+\}_{i=1}^{\infty}$ with $\tau_i \to \tilde{t}$. Consequently, we obtain:

$$v(\tau_i) - v(\tilde{t}) \geq z(\tau_i, \lambda) - z(\tilde{t}, \lambda) \tag{A5}$$

This implies:

$$D^+ v(\tilde{t}) \geq \dot{z}(\tilde{t}, \lambda) = f(\tilde{t}, z(\tilde{t}, \lambda)) + \lambda \tag{A6}$$

We distinguish the following cases:

(i) If the mapping $f(t,\cdot)$ is non-decreasing on $J \subseteq \Re$, then inequality $v(\widetilde{t}) \leq z(\widetilde{t},\lambda)$ implies $f(\widetilde{t},v(\widetilde{t})) \leq f(\widetilde{t},z(\widetilde{t},\lambda))$. The latter inequality combined with (A6) implies $D^+v(\widetilde{t}) > f(\widetilde{t},v(\widetilde{t}))$, which contradicts (2.12).

(ii) If there exists a continuous function $\phi:[t_0,T] \to \Re$ such that $f(t,w) \leq \phi(t)$, for all $(t,w) \in [t_0,T] \times J$, then we may define the lower semi-continuous function $\widetilde{v}(t) = v(t) - \int_{t_0}^{t} \phi(s)ds$. This function satisfies the following differential inequality:

$$D^+\widetilde{v}(t) \leq D^+v(t) - \phi(t) \leq f(t,v(t)) - \phi(t) \leq 0$$

Consequently, by virtue of Lemma 6.3 in [3], $\widetilde{v}(t)$ is non-increasing. This implies that $\widetilde{v}(t+h) \leq \widetilde{v}(t)$, for all $h \geq 0$. Moreover, lower semi-continuity of $\widetilde{v}(t)$ implies that for every $\varepsilon > 0$ the inequality $\widetilde{v}(t+h) \geq \widetilde{v}(t) - \varepsilon$ for sufficiently small $h \geq 0$. It follows that $\widetilde{v}(t)$ is a right-continuous function on $[t_0,t_1]$. By virtue of right-continuity and definition (A3), we must also have $v(\widetilde{t}) \geq z(\widetilde{t},\lambda)$. Thus we must have $v(\widetilde{t}) = z(\widetilde{t},\lambda)$ and in this case by virtue of (8) we obtain $D^+v(\widetilde{t}) > f(\widetilde{t},v(\widetilde{t}))$, which contradicts (2.12).

**Fact II:** $v(t) \leq w(t)$, for all $t \in [t_0,t_1)$.

Again, this claim may be shown by contradiction. Suppose that there exists $a \in (t_0,t_1)$ with $v(a) > w(a)$. Let $\varepsilon = \frac{1}{2}(v(a) - w(a)) > 0$. Furthermore, let $\lambda > 0$ be selected in such a way that (A2) is satisfied with this particular selection of $\varepsilon > 0$. Then we obtain:

$$v(a) = v(a) - w(a) + w(a) = 2\varepsilon + w(a) - z(a,\lambda) + z(a,\lambda) > \varepsilon + z(a,\lambda)$$

which contradicts Fact I. The proof is complete. ◁